# A comparison of the Benjamini-Hochberg procedure with some Bayesian rules for multiple testing[*]

## Małgorzata Bogdan[†1], Jayanta K. Ghosh[2] and Surya T. Tokdar[3]

*Wrocław University of Technology, Purdue University, Indian Statistical Institute and Carnegie Mellon University*

**Abstract:** In the spirit of modeling inference for microarrays as multiple testing for sparse mixtures, we present a similar approach to a simplified version of quantitative trait loci (QTL) mapping. Unlike in case of microarrays, where the number of tests usually reaches tens of thousands, the number of tests performed in scans for QTL usually does not exceed several hundreds. However, in typical cases, the sparsity $p$ of significant alternatives for QTL mapping is in the same range as for microarrays. For methodological interest, as well as some related applications, we also consider non-sparse mixtures. Using simulations as well as theoretical observations we study false discovery rate (FDR), power and misclassification probability for the Benjamini-Hochberg (BH) procedure and its modifications, as well as for various parametric and nonparametric Bayes and Parametric Empirical Bayes procedures. Our results confirm the observation of Genovese and Wasserman (2002) that for small p the misclassification error of BH is close to optimal in the sense of attaining the Bayes oracle. This property is shared by some of the considered Bayes testing rules, which in general perform better than BH for large or moderate $p$'s.

## 1. Introduction

Multiple tests have received considerable attention recently because of application to microarrays, where one simultaneously tests a few thousands ($m$) of null hypotheses with only a small proportion ($p$) of signals, i.e., possibly significant alternatives. Some recent references are Benjamini and Hochberg [1], Efron et al. [8], Efron and Tibshirani [7], Storey et al. [32], Genovese and Wasserman [13], Müller et al. [19], Sarkar [24] or Scott and Berger [25]. If one increases $m$ further, say $m = 10^6$, one would move from microarrays to problems of homeland security, see for example Donoho and Jin [6].

We wish to consider a still different scale, namely $m$ in the range of a few hundreds, which is relevant for quantitative trait loci (QTL) mapping. In this setup we explore and compare different multiple testing rules, ranging from the Benjamini

[*]We wish Professor P.K. Sen many more years of nonparametrics and bioinformatics.
[†]Supported by Grant 1 P03A 01430 of the Polish Ministry of Science and Higher Education
[1]Institute of Mathematics and Computer Science, Wrocław University of Technology, 50-370 Wrocław, Poland, e-mail: Malgorzata.Bogdan@pwr.wroc.pl
[2]Department of Statistics, Purdue University, West Lafayette, IN 47907-2067, USA and Indian Statistical Institute, Kolkata 700035, India, e-mail: ghosh@stat.purdue.edu
[3]Department of Statistics, Carnegie Mellon University, Pittsburgh, PA 15213, USA, e-mail: stokdar@stat.cmu.edu

and Hochberg [1] procedure [BH], Parametric Empirical Bayes [PEB] procedures, to the fully Bayes rule of Scott and Berger [25]. Also, included is a Bayesian nonparametric analysis based on Dirichlet mixtures as well as a novel application of a nonparametric algorithm for mixture estimation, due to Newton [21]. Our study is based on simulations as well as some theoretical observations. The magnitude of signals used in our simulation study is chosen according to the suggestions included in Donoho and Jin [6], so as to fulfill the condition of detectability in very sparse mixtures. For each of the considered multiple testing procedures we study its power (expected value of percentage of correctly identified alternative hypotheses), false discovery rate (FDR) and misclassification error and compare them with properties of a Bayesian oracle. We pay special attention to BH since one motivation for our study was to see if we can come up with a better Bayes or PEB rule.

Our results confirm the observation of Genovese and Wasserman [12] that for very small values of $p$'s (for $m = 200$, $p < 0.05$) the misclassification error of BH is close to optimal in the sense of attaining a Bayesian oracle. This property is shared by some of the considered Bayes testing rules, which perform better than BH for larger $p$. Moreover, in Section 3 we demonstrate that controlling positive false discovery rate (pFDR) is equivalent to controlling Bayes risk with the loss function depending only on $\alpha$ and thus, somewhat unexpectedly, the rules to control FDR or pFDR have a strong Bayesian flavor.

While our results provide some insight on QTL studies, much further work is needed to make our results directly applicable to actual QTL mapping. Our modeling is similar to that of microarrays, whereas the QTL designs require more complex linear modeling than for microarrays. The related multiple testing problem, which arises when there are many predictors (markers) to choose from, was first addressed in Bogdan et al. [4], where a suitable modification of BIC, namely mBIC, is proposed. We believe that our current research throws some light on how mBIC can be further improved by implementing a less conservative multiple testing adjustment.

The outline of the paper is as follows. In Section 2 we introduce our models and explain how some of them are related to QTL mapping. In Section 3 we discuss different notions of error in multiple testing as well as the relationship between FDR controlling rules and Bayesian testing. The *procedures considered in our study are described in Section 4, except for Bonferroni, which is described in Section 3*. The results of simulations are given in Section 5. Section 6 contains some illustrations of the problem of nonidentifiability of parameters in the mixture model and justification for using the informative prior distribution on $p$. Section 7 contains our main conclusions. Some theoretical results on the performance of the parametric Bayes procedure and the nonparametric Bayes procedure based on Dirichlet mixtures are given in the Appendix.

## 2. Models and implications for QTL mapping

We consider a multiple testing problem, when the number of tests $m$ is in the range of a few hundreds. Such values of $m$ are of importance in QTL mapping and they have a methodological interest in that the asymptotic results of Genovese and Wasserman [13], Donoho and Jin [6] or Meinshausen and Rice [18] do not yet apply.

We use the parametric model proposed in Scott and Berger [25]. Thus we consider $m$ test statistics $X_1, \ldots, X_m$ and assume that $X_i$ has either the null distribution $N(0, \sigma^2)$ or the non-null distribution $N(\mu_i, \sigma^2)$, where $\mu_i \neq 0$ represents some signal (e.g. a QTL close to the $i$-th marker). The signal $\mu_i$ is taken to be random, distributed as $N(0, \tau^2)$. Hence the non null distribution of $X_i$ is $N(0, \sigma^2 + \tau^2)$. We also define a random indicator variable $\gamma_i$, which is equal to 1 if $X_i$ is generated by the non-null distribution (i.e. it represents the signal) or 0 in the other case. If $p = P(\gamma_i = 1)$, then the marginal distribution of $X_i$ is the scale mixture of normals, namely,

$$(2.1) \qquad X_i \sim (1-p)N(0, \sigma^2) + pN(0, \sigma^2 + \tau^2).$$

Moreover, we assume that $(X_i, \gamma_i)$, $1 \leq i \leq m$, are i.i.d. random vectors. We will consider both sparse mixtures, with $p \leq 0.2$, and non-sparse mixtures, with a relatively large $p$. Usually we assume that $p$ and $\tau$ are not known, while $\sigma$ can be known or unknown, depending on the application.

For each $i$ we test whether $X_i$ has a null or non null distribution, i.e.

$$(2.2) \qquad H_{0i} : \gamma_i = 0 \text{ vs } H_{Ai} : \gamma_i = 1.$$

A major potential application of our model is QTL mapping. Our modeling takes into account the possibility that apart from QTL the trait can be influenced by a large number of polygenes, i.e. genes with very small effects, distributed over the entire genome. If our main interest is in identifying markers linked to QTL we consider a sparse mixture (2.1), where $p$ is small and $N(0, \sigma^2)$ represents the distribution of the sum of polygenic and random (environmental) effects. In this context $\sigma$ is usually unknown. The second component in the mixture, namely $N(0, \sigma^2 + \tau^2)$, represents the distribution of the QTL effect, $\mu_i$, and the sum of polygenic and random effects. Following the majority of Bayesian papers related to QTL mapping (see e.g. Yi [33]) we use $N(0, \tau^2)$ to model the distribution of $\mu_i$. Thus our model assumes that the probabilities of a positive and a negative QTL effect are the same and is suitable in the situation when the analyzed trait is not the subject of a strong selection. Note that under this scenario detecting QTL is particularly difficult. Another plausible distribution for $|\mu_i|$ is the gamma distribution (see e.g. Otto and Jones [22]). A completely robust alternative is to model $\mu_i$'s with a non-parametric distribution $P$ and put a further prior $P \sim$ Dirichlet, which leads to Dirichlet location mixture distribution for $X_i$'s. We investigate this approach and propose an alternative nonparametric inference based on Newton [21].

If our main interest is in both QTL and polygenic effects, the null component, $N(0, \sigma^2)$, represents the distribution of random effects, and $N(0, \tau^2)$, represents the distribution of effects due to QTL and polygenes. In this setting $p$ need not be small and $\sigma^2$ may be assumed known, since we can precisely estimate it through replications.

**Remark 1.** The number of strong QTL, which are significantly different from the background of polygenes, is usually small. In this case only relatively large QTL effects, $|\mu_i| \geq \sigma\sqrt{2 \log m}$, may be identified, since extreme values of the "null" component of the mixture are approximately equal to $\sigma\sqrt{2 \log m}$. In order that such signals are generated by the non null component $\tau^2$ should be comparable or larger than $2\sigma^2 \log m$.

**Remark 2.** The assumption of the independence of $X_i$ can be used when markers are distant from each other. When markers are close to each other, the corresponding test statistics might be strongly correlated. However, the results reported in this paper demonstrate some general properties of the multiple testing procedures and show the directions in construction of related methods for detection of linked QTL.

## 3. Different notions of error in multiple testing

Consider the problem of testing of $m$ hypothesis $H_{01}, \ldots, H_{0m}$, specified in (2.2). For each individual test two types of error can occur: the null hypothesis can be rejected even though it is true (type I error) or be not rejected when it is false (type II error). Following the notation of Benjamini and Hochberg [1], Table 1 defines variables describing counts of possible outcomes of a multiple testing procedure.

The main focus of classical statistics is on tests minimizing the probability of the type II error (or maximizing the power), while controlling the probability of the type I error at a given significance level $\alpha$. The natural extension of the type I error to the situation of testing $m$ hypotheses is the family-wise error rate, $FWER = P(V > 0)$. Additionally, the notion of power can be naturally extended to the multiple testing as $E(\frac{S}{m_1}|m_1 > 0)$. Here, as well as in the next part of the paper, $E$ is used to denote the frequentist expectation (i.e. conditional on the vector parameters of the model (2.1)). The classical approach to the multiple testing problem relies on constructing procedures maximizing the power while controlling FWER at a given level (see e.g. Holm [16]).

In the situation when $m$ is large, procedures controlling FWER are usually very conservative. Note that in many practical applications one would often accept false discoveries as long as they consist only a small proportion of all discoveries. Going along these practical expectations Seeger [26] elaborated on the idea of Eklund (unpublished seminar papers) and discussed a stepwise multiple testing procedure aimed at controlling the proportion of false discoveries among all discoveries. The same stepwise multiple testing procedure has been later discovered by Simes [27], who proved that it controls FWER in a weak sense (when all hypothesis are true). The notion of proportion of false discoveries appeared again in a paper by Sorić [28]. Following this paper, Benjamini and Hochberg [1] formally defined the false discovery rate as $FDR = E(\frac{V}{R})$, where $\frac{V}{R} = 0$ if $R = 0$. Benjamini and Hochberg also proved that the multiple testing procedure of Seeger and Simes controls FDR at a desired level when the test statistics are independent. Following Benjamini and Hochberg [1] this procedure gained a large popularity and is currently known as the Benjamini and Hochberg (BH) procedure.

Let $P_{(1)} \leq P_{(2)} \ldots \leq P_{(m)}$ be the ordered p-values of $m$ tests. Let

$$(3.1) \qquad k = \max\left\{i: \ P_{(i)} \leq \frac{i\alpha}{m}\right\}.$$

BH rejects all hypotheses for which the corresponding p-values are smaller than $P_{(k)}$. In Benjamini and Yekutieli [3] and Sarkar [23] it is proved that BH controls FDR also under certain forms of positive dependence between test statistics. Following Benjamini and Hochberg [1] many other criteria and procedures which allow for controlling a number or proportion of false discoveries, were developed (see e.g. Lehmann and Romano [17], Sarkar [24], Storey [31] and references given there) but BH still remains one of the most popular methods of multiple testing.

TABLE 1
*Counts of possible outcomes of m hypothesis tests*

|  | Accept null | Reject null | Total |
|---|---|---|---|
| Null true | U | V | $m_0$ |
| Alternative true | T | S | $m_1$ |
|  | W | R | $m$ |

Multiple testing problems can be approached also from the point of view of decision theory. Depending on the specifics of the problem, different loss functions can be assigned to the two types of errors and the procedure minimizing the related risk can be constructed. The corresponding procedures in the framework of Bayesian decision theory were discussed e.g. in Müller et al. [19], Müller et al. [20] or Scott and Berger [25]. Further references to Bayesian multiple testing procedures as well as a novel Bayesian stepwise multiple testing procedure can be found in Chen and Sarkar [5].

To point at some similarities between controlling FDR and Bayesian approach to multiple testing we now briefly discuss the positive false discovery rate,

$$pFDR = E\left(\frac{V}{R}|R>0\right) = \frac{FDR}{P(R>0)},$$

defined in Storey [30] and Bayesian false discovery rate

$$BFDR = P(H_0 \text{ is true}|H_0 \text{ is rejected}),$$

defined in Efron and Tibshirani [7].

Theorem 1 of Storey [30] states that in case when individual test statistics are generated by the two-component mixture model, like in our setting, $pFDR = BFDR$. It is also pointed out that there are situations in which BFDR can not be controlled. An obvious example is when $p = 0$, since then $BFDR = pFDR = 1$. It is however easy to show that in our testing problem (2.2) $BFDR$ can be controlled at any given level $\alpha$ if $p > 0$. The corresponding threshold for the absolute value of the test statistic $X_i$ is given by the formula

$$(3.2) \qquad c_{fdr} = \inf\left\{x > 0 \ : \ \frac{(1-p)(1-\Phi_0(x))}{1-F(x)} < \alpha\right\},$$

where $\Phi_0$ and $F$ are cdfs of $N(0, \sigma^2)$ and the mixture distribution (2.1), respectively.

**Remark 3.** The difference between FDR and BFDR may be relatively large for small $p$ and a small deviation between the null and alternative distribution (i.e. small power). However, in typical QTL or microarray experiments, where $m$ is large and some rejections typically occur, the difference between BFDR and FDR is usually very small. Based on the asymptotic approximation of FDR by BFDR, Genovese and Wasserman [13] call (3.2) an *oracle threshold* to control FDR.

**Remark 4.** Theorem 5.1 of Benjamini and Yekutieli [3] states that if the test statistics are continuous and independent then FDR of BH is equal to $\alpha m_0/m$. Thus FDR of BH is close to $\alpha$ only when $m_0$ is close to $m$ and converges to 0 when $m_0 \to 0$. When $m_0$ is known one can easily modify BH to control FDR at the level $\alpha$ by replacing $k$ (see 3.1) with

$$(3.3) \qquad k1 = \max\left\{i \ : \ P_{(i)} \leq \frac{i\alpha}{m_0}\right\}.$$

In Benjamini and Hochberg [2] a graphical method to estimate $m_0$ is proposed and the formula (3.3) is used to construct an adaptive version of BH.

Under the mixture model (2.1) the expectation of $m_0$ is equal to $m(1-p)$ and a corresponding modified version of BH can be obtained by replacing $k1$ with

$$(3.4) \qquad k2 = \max\left\{i \ : \ P_{(i)} \leq \frac{i\alpha}{m(1-p)}\right\}.$$

TABLE 2
*Matrix of losses*

|              | Accept $H_{0i}$ | Reject $H_{0i}$ |
|--------------|-----------------|-----------------|
| $H_{0i}$ true | 0              | $\delta_0$      |
| $H_{Ai}$ true | $\delta_A$     | 0               |

It is easy to prove that this version of BH also has FDR equal to $\alpha$. Moreover, in Efron and Tibshirani [7] it is noticed that the modified BH (3.4) is equivalent to the BFDR controlling rule (3.2), with the cdf of the mixture distribution estimated by the empirical distribution function. In many consecutive papers (see e.g. Efron et al. [8], Efron and Tibshirani [7], Storey [29] and Genovese and Wasserman [13]) different nonparametric methods of the estimation of $(1-p)$ and $F(x)$ were considered, leading to FDR controlling rules which are more liberal than BH.

Let us now consider the multiple testing problem from the perspective of decision theory. Table 2 defines the specific matrix of losses for making the wrong decision.

Let us denote by $t_1$ and $t_2$ the probability of type I and type II errors of a single test. The Bayes risk related to the above matrix of losses is given by the following equation

$$BR_{\delta_0,\delta_A} = \delta_0(1-p)t_1 + \delta_A p t_2. \tag{3.5}$$

The Bayes rule, i.e., the test which minimizes this risk, rejects the null hypothesis if

$$\frac{f_A(X_i)}{f_0(X_i)} > \frac{(1-p)\delta_0}{p\delta_A}, \tag{3.6}$$

where $f_0$ and $f_A$ are the densities of $X_i$ under $H_0$ and $H_A$, or equivalently if

$$p_i = P(H_{Ai}|X_i) > \frac{\delta_0}{\delta_0 + \delta_A}. \tag{3.7}$$

We call this test a Bayes oracle and compare other tests to this oracle.

Let us observe that

$$BFDR = \frac{(1-p)t_1}{(1-p)t_1 + p(1-t_2)}.$$

Thus

$$BFDR < \alpha \text{ iff } (1-\alpha)(1-p)t_1 + \alpha p t_2 < \alpha p \tag{3.8}$$

and controlling BFDR controls the Bayes risk with a loss $\delta_0 = 1 - \alpha$ and $\delta_A = \alpha$. The classical flavor of BFDR is however strongly reflected in assigning much larger loss to the type I error than to the type II error.

The accuracy of the multiple testing procedure can be judged by its misclassification probability, $MP = \frac{E(V+T)}{m}$. Note that $MP = BR_{1,1}$, where $BR_{1,1}$ is the Bayes risk corresponding to 0-1 loss. In our parametric setting (2.1) the Bayes oracle minimizing $BR_{1,1}$ rejects the hypothesis $H_{0i}$ if

$$X_i^2 > \frac{2(\sigma^2 + \tau^2)\sigma^2}{\tau^2}\left(\frac{1}{2}\log\left(\frac{\sigma^2 + \tau^2}{\sigma^2}\right) + \log\frac{1-p}{p}\right). \tag{3.9}$$

In Figure 1 we compare the Bayes oracle (3.9) to BH and the standard testing procedure based on the Bonferroni correction. The significance level for each individual test in Bonferroni procedure is equal to $\alpha/m$. For this presentation as well as for simulations reported in Section 5 we use $m = 200$, $\sigma = 1$ and $\tau = \sqrt{2 * \log(200)} \approx 3.26$. For BH and Bonferroni procedures $\alpha = 0.05$. Apart from the standard BH we use its modified version, with the cutoff for p-values $k2$ specified by (3.4). The reported characteristics for the Bonferroni correction and Bayes oracle were obtained theoretically, while the characteristics of BH were computed using computer simulations, based on 10000 replicates.

Figure 1 demonstrates that, as expected, the modified version of BH keeps FDR exactly at the level 0.05, while FDR of the original BH decreases linearly with $p$. Comparison of 1(a) and 1(b) shows that the difference between BFDR and FDR is substantial when $p < 0.03$. In particular, neither versions of BH controls BFDR in this range of $p$. This seems due to the fact that for very small $p$ the threshold based on the empirical mixing distribution is substantially more liberal than the one provided by (3.2). Both versions of BH take an intermediate position between the Bonferroni procedure, which is most conservative, and the most liberal Bayes oracle. Figure 1 demonstrates that the most powerful Bayes oracle has also the largest FDR. However, as expected, type I and type II errors balance in such a way that the misclassification probability ($MP$) of the Bayes oracle is smaller than that of any other method. Interestingly, the modified version of BH performs very well in terms of $MP$ over the entire range of $p$. When $p$ is very small also the original BH has a very low MP, which for $p = 0.015$ is very close to the optimal value provided by the Bayes oracle.

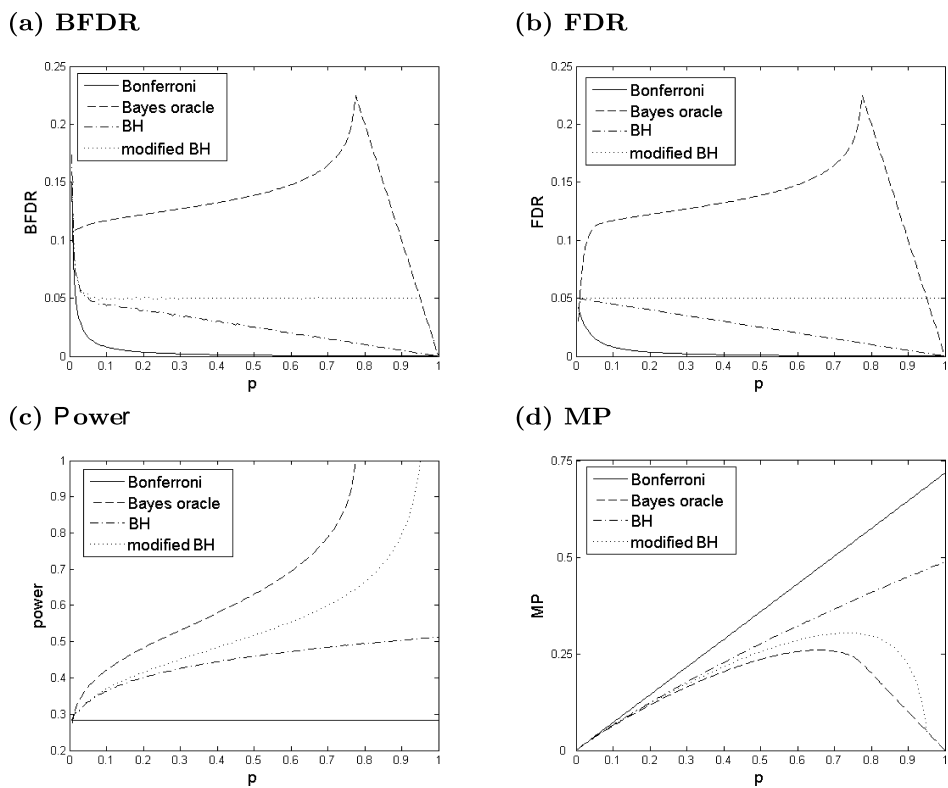

FIG 1. *Characteristics of multiple testing procedures.*

## 4. Bayes, parametric empirical Bayes and modified BH procedures

**PEB procedures:** The natural way of applying the Bayes classifier (3.9) in the situation when parameters of (2.1) are unknown is to use their consistent estimates and plug them into (3.9). In particular, maximum likelihood estimators (MLE) could be considered.

Let

$$L(X_1, \ldots, X_m | p, \tau, \sigma) = \prod_{i=1}^{m} (p f_A(X_i) + (1-p) f_0(X_i)).$$

We estimate our model parameters in two steps. First we fix $p$ and estimate $\tau(p)$ and $\sigma(p)$ using the EM algorithm. In the second step we estimate $p$ by maximizing $L(X_1, \ldots, X_m | p, \hat{\tau}(p), \hat{\sigma}(p))$ using numerical methods. We plug our estimated parameters into (3.9) and *denote the resulting Parametric Empirical Bayes Classifier as PEB1.* As reported in Section 5, PEB1 performs very well for moderate values of $p$. However, when $p$ is very small PEB1 has large FDR and MP. This behavior is related to the problems with identifiability of parameters of mixture distributions discussed in Section 6. Since our main interest is in sparse mixtures we consider the following modification of PEB1. Firstly we stabilize the performance of MLE by supplying the information included in the data with the prior information on $p$. Using a subjective, informative, prior on $p$ is also strongly recommended in Scott and Berger [25], where the following prior density is proposed;

$$(4.1) \qquad f(p) = \beta(1-p)^{\beta-1}.$$

In simulations reported in Scott and Berger [25], the parameter $\beta$ is set to be equal to 11, so the corresponding prior on $p$ has its median close to 0.07. To adjust to the sparsity typical for QTL mapping experiments we slightly shift this prior towards 0 and choose $\beta$ so that the prior median is 0.03, which for $m = 200$ corresponds to 6 signals on average. In our simulations $\beta = 22.76$. The results aren't sensitive to small changes in $\beta$.

The estimate of $p$ is obtained by maximizing

$$(4.2) \qquad \log L(X_1, \ldots, X_m | p, \hat{\tau}(p), \hat{\sigma}(p)) - (\beta - 1) \log(1-p)$$

and can be interpreted as a mode of the "posterior" density of $p$.

The second modification relies on replacing the maximum likelihood estimates of $\tau(p)$ and $\sigma(p)$ with the moment estimates based on the fourth and the second moment of the mixture distribution. We observe that using the fourth moment makes our procedure sensitive to the change in the tail of the mixture distribution and hence yields good results in a very sparse mixture case. *The resulting Parametric Empirical Bayes Classifier is denoted by PEB2.*

When $\sigma$ is known PEB1 and PEB2 are constructed accordingly. For PEB2 $\tau(p)$ is estimated using the fourth moment of the mixture distribution.

**Modified BH:** We use estimates of $p$ and $\sigma$ computed by PEB methods to construct modified versions of BH, with the threshold specified by (3.4). The version based on MLE is denoted by BH1 and the sparse mixture version, based on the estimates derived by maximizing (4.2), by BH2.

**Full Bayes approach:** We use the framework of Scott and Berger [25] and construct the full Bayes procedure minimizing the posterior Bayes risk corresponding to the 0-1 loss (1 for making type I or type II error). We use noninformative priors for $\tau^2$ and $\sigma^2$, suggested in Scott and Berger [25], with densities

$$(4.3) \qquad \pi_{SB}(\sigma^2) = \frac{1}{\sigma^2} \text{ and } \pi_{SB}(\tau^2|\sigma^2) = \frac{\sigma^2}{(\sigma^2+\tau^2)^2}.$$

When $\sigma$ is known only the prior for $\tau$ is used. The prior on $p$ is the same as the one used by PEB methods, namely (4.1), with $\beta = 22.76$. To compute the posterior probability of $H_{0i}$ (see formula (9) of Scott and Berger [25]) Markov chain Monte Carlo (MCMC) is applied, which according to our simulations is more stable than the importance sampling suggested in Scott and Berger [25]. The hypothesis $H_{0i}$ is rejected if

$$(4.4) \qquad P(H_{0i}|X_1,\ldots,X_m) < 0.5$$

and the resulting multiple testing procedure is denoted by SB.

**Remark 5.** Note that minimizing the posterior Bayes risk is conceptually different from minimizing the risk $BR_{1,1}$ (see (3.5)). $BR_{1,1}$ is conditional on the vector of parameters of the mixture model (2.1), while the posterior Bayes risk is conditional on the data and depends on the prior. However, Theorem 8.2 in the Appendix states that if the parameter space (i.e. $p \in (0,1)$, $\sigma > 0$ and $\tau > 0$) and $m$ increased then the misclassification probability of SB converges to the optimal value provided by the oracle (3.7). This result is a consequence of Theorem 8.1 on posterior consistency under the considered mixture model. Obviously, for each fixed $m$, the difference between misclassification probability of SB, and the oracle depends on the accuracy of prior assumptions and due to the choice of the prior on $p$ we expect SB to resemble the oracle when the data are generated by the sparse mixture.

**Dirichlet mixtures:** The procedures presented so far are based on the assumption that the distribution of the signals (of $\mu_i$'s given $\gamma_i = 1$) is completely known up to finitely many parameters. In practice, however, a lot less is known about the signals. A realistic model for such a situation is to consider : $\mu_i|(\gamma_i = 1) \sim P^{sig}$ for some unknown probability measure $P^{sig}$ with $P^{sig}(\{0\}) = 0$, which doesn't need to be restricted to any parametric family. In this case, $X_i$'s arise as independent observations from the mixture density $f(x) = \int \phi_\sigma(x-\mu)dP(\mu)$, where $P$ is a probability measure that puts some positive mass at the point 0 and distributes the remaining mass $p$ according to $P^{sig}$. A Bayesian analysis of this model is possible by using a prior distribution on the space of such probability measures $P$. A suitable candidate is a Dirichlet process prior. Below we introduce a new procedure based on this model and prior.

Let $\pi_{SB}(\tau^2, \sigma^2) = (\sigma^2+\tau^2)^{-2}$ denote the joint prior distribution on $(\tau^2, \sigma^2)$ recommended in Scott and Berger [25]. We assume that

$$\begin{aligned}
X_i|\mu_i, P, \sigma^2, \tau^2, p_0, c &\sim N(\mu_i, \sigma^2), \\
\mu_i|P, \sigma^2, \tau^2, p_0, c &\sim P, \\
P|\sigma^2, \tau^2, p_0, c &\sim Dir(c, (1-p_0)\delta_{\{0\}} + p_0 N(0, \tau^2)), \\
(p_0, c, \tau^2, \sigma^2) &\sim Beta(1, 22.76) \times Gamma(1,1) \times \pi_{SB}(\tau^2, \sigma^2),
\end{aligned}$$

where $Dir(c, P_0)$ denotes the Dirichlet process prior (see Ferguson [11]) with precision constant $c > 0$ and base measure $P_0$ - a probability measure on the real line. Our choice of the base measure, namely $P_0 = (1-p_0)\delta_0 + p_0 N(0, \tau^2)$, ensures that a random $P \sim Dir(c, P_0)$ almost surely puts some positive mass $1-p$ on 0 and distributes the remaining positive mass $p$ on the real line according to some probability measure $P^{sig}$, which is singular to $\delta_{\{0\}}$. Therefore, without any ambiguity, we can import our familiar signal indicators $\gamma_i$ into this model by defining $\gamma_i = I(\mu_i \neq 0)$.

Note that the priors for $p_0$ and $(\tau^2, \sigma^2)$ match the priors chosen for these parameters in the model proposed in Scott and Berger [25] and presented in the previous section. The precision constant $c$ is modeled with a $Gamma(1,1)$ prior, which is quite diffuse with a mean equal to 1 – a conventional choice of this parameter.

Toward implementation of this model, we first integrate out $P$ from the hierarchical structure by using the Pólya urn representation of a Dirichlet process. Although our specification includes an improper prior on $\sigma^2$, it turns out that the posterior distribution of $(\mu_1, \ldots, \mu_m, \sigma^2, \tau^2, p_0)$ is indeed a proper distribution; see Theorem A.3 in Appendix. This allows us to obtain an MCMC sample of observations $(\mu_1^{(l)}, \ldots, \mu_m^{(l)})$, $l = 1, \ldots, L$, from the posterior distribution of $\mu_i$'s given the data. We use the algorithm described in Escobar and West [10] with suitable adaptations to our model. Our model differs from the one considered in Escobar and West [10] in two aspects: 1. we consider only a location mixture of normals and 2. our base measure has a point mass at $\{0\}$. The adaptations, however, are not complicated and we omit further details.

With the sample of $\mu_i$'s collected from the posterior we calculate $P(\gamma_i = 0|X_1, \ldots, X_n) \approx \frac{1}{L}\sum_l I(\mu_i^{(l)} = 0)$. As before, we reject $H_{0i}$ if this estimate is smaller than 0.5. We denote this multiple testing procedure by DPP.

**Approximate nonparametric Bayes procedure based on Newton's algorithm, NPBN:** A somewhat related procedure can be obtained by combining the above nonparametric model with Newton's algorithm (see Newton [21]), which produces an easy to compute, recursive estimate of the distribution $P$. In particular, we start with an initial guess of $P$ given by $P_0 = (1-p_0)\delta_{\{0\}} + p_0 N(0, \tau^2)$ and then recursively update this guess as

$$P_i(d\mu) = (1-w_i)P_{i-1}(d\mu) + w_i \frac{\phi_\sigma(x_i - \mu)P_{i-1}(d\mu)}{\int \phi_\sigma(x_i - \nu)P_{i-1}(d\nu)},$$

where $w_i \in (0,1)$ are prefixed weights (we take $w_i = (i+1)^{-1}$). We take the final update $P_m$ as the estimate of $P$. Note that the estimate $P_m$, too, puts some positive mass $1 - p_m$ at 0 and distributes the rest according to some density $f_m$ on the real line. Testing is then performed by mimicking the Oracle rule and replacing $P$ with the estimate $P_m$: we reject $H_{0i}$ if $(1/p_m - 1)\phi_\sigma(x_i)/\int \phi_\sigma(x_i - \mu)f_m(\mu)d\mu < 1$. We call this procedure Nonparametric Bayesian Procedure based on Newton's algorithm (NPBN).

For every $i$, if one models $x_i \sim \int \phi_\sigma(x - \mu)P(d\mu)$ with $P \sim Dir(1, P_{i-1})$ then the posterior expectation of $P$ given the singleton sample $\{x_i\}$ equals $P_i$. In spite of this resemblance, NPBN should not be taken as an approximation to DPP. The former, however, has its own set of advantages.

The biggest advantage of using NPBN is that it produces extremely fast computation while using a nonparametric model. The reason for its speed stems from the *one pass routine* employed by the algorithm.

The output of the NPBN procedure depends on the order in which data are fed to the algorithm. In our simulations we align the observations in their ascending order of magnitude. With this alignment, $P_i$'s are first trained on small observations, which are mostly noise, followed by the large ones coming mainly from the signals. However, as the later updates are rather less influential (small $w_i$), the concentration of $P_m$ on 0 would be systematically inflated. Therefore the chosen alignment would systematically result in a more conservative procedure than what a random alignment would produce. Such a conservative approach is well suited to our anticipation of a small to moderate number of signals.

While greater speed is a selling point for NPBN, it does suffer from some inflexibility in model specification. Unlike DPP, the NPBN setting does not allow a further prior specification on the parameters $p_0, \tau^2, \sigma^2$. It is hard to generalize the recursive algorithm to include unknown parameters. In the present paper we consider NPBN only when $\sigma$ is known and specify $\tau^2 = \sigma^2$, which is equal to the mean of the prior distribution on $\tau^2$ used by SB and DPP.

We choose $p_0$ as

$$p_0 = \frac{1}{m} \sum_i I\left(\frac{\phi_\sigma(x_i)}{\phi_{\sqrt{\sigma^2+\tau^2}}(x_i)} < 1\right).$$

This quantity is equal to the proportion of rejections one would make by assuming $P = \frac{1}{2}\delta_0 + \frac{1}{2}N(0, \tau^2)$. This choice calibrates $p_0$ to $\tau^2$ in a natural way - once $\tau^2$ is picked we update our noninformative choice of $p_0 = 1/2$ by using this chosen value of $\tau^2$. Our simulation study indicates that this data dependent choice of $p_0$ leads to an overall higher efficiency compared to any fixed choice of $p_0$.

## 5. Simulation results

Table 3 and Figure 2 demonstrate characteristics of SB, PEB1, PEB2, BH1, BH2 and NPBN. "Efficiency", represented in Figures 2(a) and (c), is defined as

$$E = \frac{\text{MP of the oracle}}{\text{MP of a given procedure}}.$$

We do not report the results of the original BH since the performance of BH2 is systematically better. The parameter values used in the simulations are $m = 200$, $\sigma = 1$ and $\tau = \sqrt{2 * \log(200)} \approx 3.26$. Due to the computational complexity the results for SB are based on 3000 replicates. The results of all other procedures are based on 10000 replicates. The large scale simulations were not feasible for DPP, which is not represented in Table 3 and Figure 2.

Table 3 demonstrates that for $p \leq 0.05$ PEB1 and BH1 have large MP and FDR. The properties of these rules quickly improve with increasing $p$ and for $p \approx 0.2$ MP of PEB1 is close to optimal and FDR of BH1 is close to 0.05. When $\sigma$ is known the characteristics of PEB1 and BH1 stay at the assumed level for all $p \geq 0.2$ but when $\sigma$ is unknown these rules deteriorate again when $p > 0.8$. The sparse mixture rules: SB, PEB2 and BH2, perform well for very small $p$. When $\sigma$ is known these rules retain good properties for $p \in [0, 0.6]$ but when $\sigma$ is unknown they deteriorate already at $p \approx 0.3$. Figure 2(d) demonstrates that at this point all sparse mixture rules start to loose power and become overly conservative. The reason for this behavior as well as the corresponding loss of power for PEB1 and BH1 when $p > 0.8$ is the difficulty with identifying the model parameters, discussed in detail in Section 6.

Figures 2(a) and (c) demonstrate the "efficiencies" of the sparse mixture rules in the most interesting range $p \leq 0.2$. PEB1 and BH1 are not represented since their "efficiencies" for $p < 0.03$ are below 50%. Figures 2(a) and (c) show that when $p \in [0.01, 0.03]$ BH2 is almost optimal and has the "efficiency" slightly larger than the "efficiencies" of other procedures. However, this characteristic of BH2 systematically decreases with an increase of $p$ and at $p = 0.2$ it is substantially smaller

TABLE 3
*FDR and Misclassification Probability of multiple testing procedures. BO stands for the Bayes Oracle (3.9)*

| | $\sigma$ known | | | | | | | $\sigma$ unknown | | | | |
|---|---|---|---|---|---|---|---|---|---|---|---|---|
| $p$ | BO | SB | PEB1 | PEB2 | BH1 | BH2 | NPBN | SB | PEB1 | PEB2 | BH1 | BH2 |
| | Misclassification probability in % | | | | | | | | | | | |
| 0.0 | 0 | 0.01 | 82.4 | 0.04 | 73.2 | 0.03 | 0.01 | 0.02 | 31.3 | 0.04 | 14.3 | 0.03 |
| 0.025 | 1.76 | 1.8 | 19.3 | 1.77 | 16.8 | 1.77 | 1.82 | 1.84 | 11.0 | 1.81 | 5.11 | 1.80 |
| 0.05 | 3.36 | 3.38 | 7.01 | 3.40 | 6.19 | 3.42 | 3.46 | 3.48 | 6.36 | 3.48 | 4.35 | 3.51 |
| 0.2 | 11.7 | 11.8 | 11.8 | 11.8 | 12.1 | 12.2 | 11.9 | 12.4 | 12.2 | 12.3 | 12.3 | 13.0 |
| 0.5 | 23.5 | 24.5 | 24.1 | 24.5 | 25.5 | 26.2 | 24.0 | 35.5 | 24.9 | 40.0 | 26.4 | 42.6 |
| 0.8 | 20.0 | 29.6 | 21.1 | 29.5 | 28.8 | 35.1 | 22.2 | 79.6 | 30.2 | 79.7 | 41.7 | 79.9 |
| | False Discovery Rate in % | | | | | | | | | | | |
| 0.0 | 0 | 3.1 | 93.1 | 6.0 | 79.4 | 5.2 | 2.4 | 4.4 | 49.4 | 5.5 | 42.4 | 5.6 |
| 0.025 | 9.4 | 7.8 | 31.5 | 7.2 | 20.5 | 5.0 | 5.5 | 8.6 | 28.0 | 6.3 | 18.7 | 4.1 |
| 0.05 | 11.2 | 8.7 | 17.5 | 8.0 | 8.0 | 4.9 | 6.9 | 9.2 | 19.9 | 7.0 | 11.0 | 3.9 |
| 0.2 | 12.2 | 9.2 | 12.9 | 8.4 | 5.0 | 4.7 | 8.6 | 5.5 | 13.9 | 7.6 | 6.5 | 2.9 |
| 0.5 | 13.9 | 8.0 | 14.7 | 8.0 | 5.3 | 4.1 | 11.0 | 1.1 | 14.6 | 2.0 | 6.2 | 0.7 |
| 0.8 | 20.0 | 5.0 | 17.2 | 5.1 | 6.1 | 2.4 | 14.8 | 0.0 | 13.1 | 0.0 | 4.4 | 0.0 |

### $\sigma$ known

**(a) E**

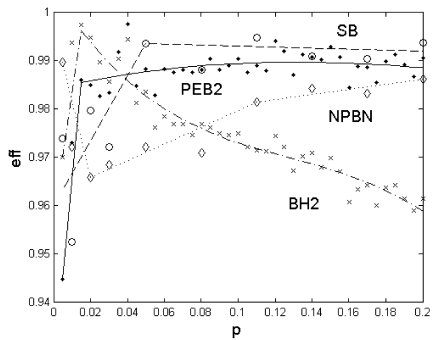

**(b) Power**

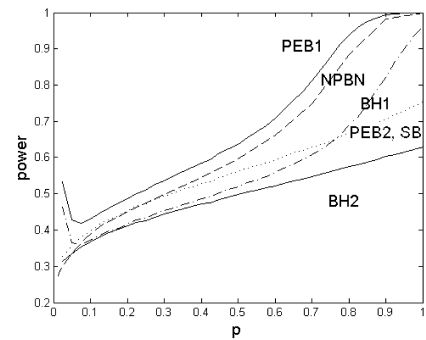

### $\sigma$ unknown

**(c) E**

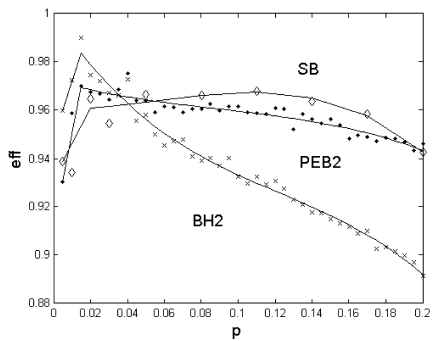

**(d) Power**

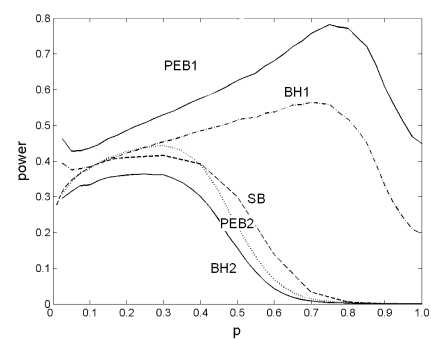

FIG 2. *Characteristics of multiple testing procedures.*

than the "efficiencies" of PEB2 and SB. The "efficiencies" of PEB2 and SB stay constant at the level close to 99% when $\sigma$ is known and only slightly decrease to 94.5% at $p = 0.2$ when $\sigma$ is unknown.

When $\sigma$ is known the nonparametric NPBN procedure performs surprisingly well over the entire range of $p$. It is slightly less efficient than SB and PEB2 when $p < 0.2$ but it retains good properties also for $p$ close to 1.

In order to get a feeling of the performance of DPP we compare it with SB on a case by case basis with the help of a few toy data sets. We generated 10 data sets, each of size 200, from the model described in Section 1, with $\sigma = 1$, $\tau = \sqrt{2 \log 200}$ and various values of $p$. Five of these data sets are represented in Figure 3. On each panel, two scatter plots of $\Pr(\gamma_i = 0 | x_1, \ldots, x_m)$ versus $x_i$ are presented. The open circles joined by the solid line correspond to DPP, whereas the filled circles joined by the dotted line correspond to SB. The left column in Figure 3 corresponds to the known $\sigma$ case – i.e., both DPP and SB are employed with $\sigma^2$ fixed at 1 and the conditional prior $\pi(\tau^2 | \sigma^2 = 1)$ is used for $\tau^2$. The right column corresponds to the unknown $\sigma$ case.

It appears that DPP and SB perform quite similarly, although the former is a little more conservative than the latter, particulary when the number of signals is very small. This is further illustrated by Table 4. The ten columns in the table represent the ten data sets, with the number of signals shown on the header row. In each cell of the body of the table, the two values give the numbers of correct and incorrect discoveries of signals made by the corresponding procedure for that particular data set. From Table 4 we also note that DPP and NPBN are quite similar except for samples with many signals (last two columns), where the prior on $p$ used by DPP was strongly inappropriate.

## 5.1. When the prior assumption is wrong

In this section we demonstrate the results of simulations illustrating the performance of our methods in the situation when the assumption that the distribution of $\mu_i$ under the alternative is normal does not hold. For this simulation we consider the case with $\sigma$ known and generate $\mu_i$'s using a symmetrized gamma distribution instead of normal distribution.

Let $g(x, r, u)$ denote the density of the gamma distribution with the shape parameter $r$ and the scale parameter $u$. The symmetrized gamma density describing the prior distribution of $\mu_i$ under the alternative is given by the equation $g_A(x) = 0.5 g(|x|, r, u)$. For the current simulation we use $r = 4$ and $u = \dfrac{2\sqrt{2 \log m}}{r\sqrt{2\pi}}$.

As demonstrated in Figure 4, the deviation from the prior assumption strongly affects the behavior of PEB1 and BH1, which are based on the maximum likelihood estimates of parameters under the wrong mixture model (2.1). In particular, for $p \in (0.3, 0.8)$ these procedures are much too liberal and do not achieve the assumed characteristics. The misclassification probability of PEB1 is close to optimal only in the range of $p \in (0.1, 0.2)$ and $p > 0.9$. Also, only in this range FDR of BH1 is close or below the assumed value of 0.05. Over the entire range of $p$ the misclassification probability of the nonparametric procedure NPBN is decisively smaller than for PEB1, which clearly demonstrates the advantage of using nonparametric methods in case when the prior distribution is not known. Surprisingly, PEB2 and BH2, which are based on the moments estimates under the wrong model and use the strongly informative prior on $p$, behave very well over the entire range of $p$. We believe that their good behavior for moderate and large $p$ is just a coincidence, resulting from the opposite influence of different types of errors of our estimation procedure. As noted before, under this particular violation of the prior assumption the methods

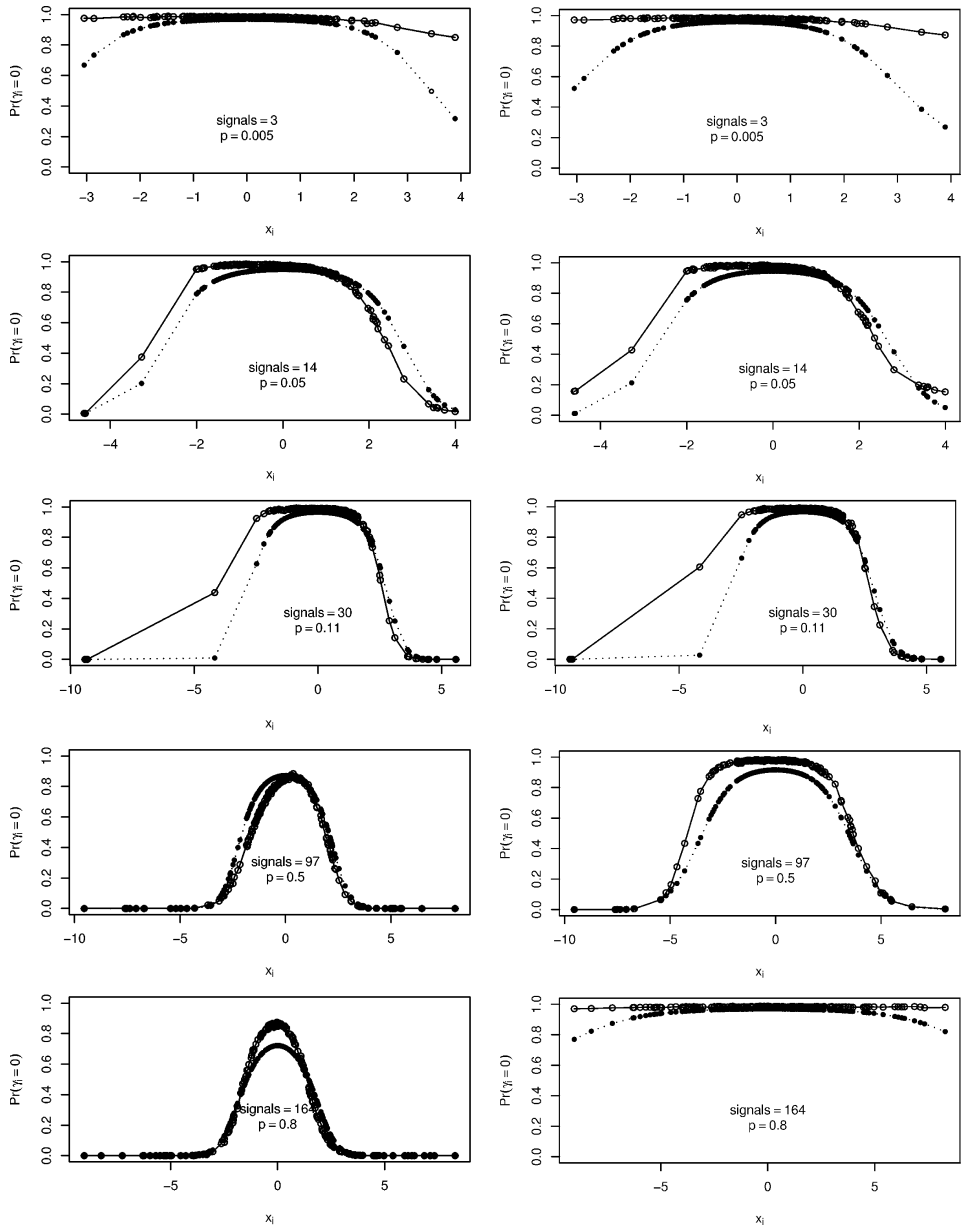

FIG 3. *A case by case comparison of DPP and SB for known $\sigma$: scatter plots of $\Pr(\gamma_i = 0|x_1,\ldots,x_m)$ vs. $x_i$. Open circles connected by solid line represent DPP, filled circles connected by dotted line represent SB. Left column depicts the situation when $\sigma$ is known, right when $\sigma$ is unknown.*

which do not use the prior on $p$ are too liberal when $p \in (0.3, 0.8)$. In case of PEB2 and BH2 this error seems to cancel the error resulting from using the informative, but not adjusted to this range, prior on $p$. However, other simulations, not reported in this paper, suggest that the good behavior of PEB2 and BH2 for $p < 0.2$ is a quite general rule, working under a wide set of different, also asymmetric, prior distributions on $\mu_i$. A theoretical explanation of this phenomenon still needs to be

Table 4
*A comparison of the numbers of correct and false discoveries of DPP, SB and NPBN*

| signals | 3 | 5 | 14 | 12 | 19 | 30 | 97 | 100 | 164 | 146 |
|---|---|---|---|---|---|---|---|---|---|---|
| | | | | | $\sigma$ known | | | | | |
| SB | 2 0 | 2 1 | 9 1 | 6 0 | 5 2 | 14 0 | 51 3 | 60 7 | 110 10 | 95 6 |
| DPP | 0 0 | 0 0 | 9 3 | 6 0 | 7 2 | 14 0 | 62 7 | 58 7 | 109 9 | 89 2 |
| NPBN | 1 0 | 1 0 | 8 1 | 4 0 | 3 0 | 15 1 | 59 7 | 64 7 | 164 36 | 111 18 |
| | | | | | $\sigma$ unknown | | | | | |
| SB | 2 0 | 3 2 | 9 1 | 8 3 | 2 0 | 14 0 | 29 0 | 56 5 | 0 0 | 0 0 |
| DPP | 0 0 | 0 0 | 9 2 | 7 0 | 0 0 | 13 0 | 23 0 | 39 0 | 0 0 | 0 0 |

**(a) MP**  **(b) FDR**

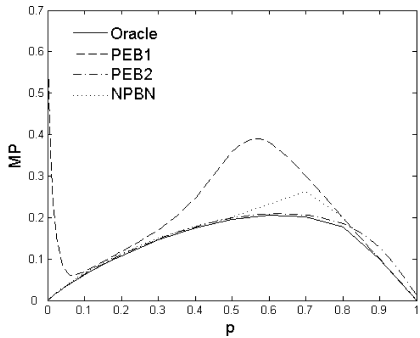 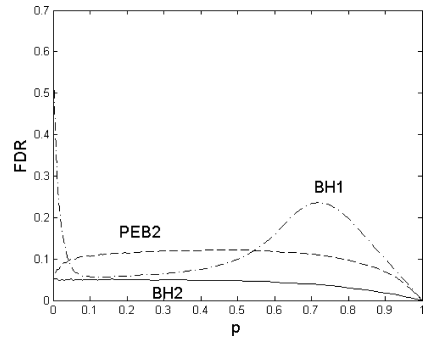

FIG 4. *Characteristics of multiple testing procedures when the assumption on the prior distribution of $\mu_i$ is violated.*

worked out.

## 6. Problems with identifiability of model parameters

An early treatment of the problem caused by lack of identifiability of mixtures can be found in Ghosh and Sen [15]. Another recent reference is Elmore et al. [9]. In this section we illustrate this problem by a numerical study on the Kullback-Leibler distance between different mixture densities and the resulting behavior of the maximum likelihood method.

Consider the problem of a choice between two competing probability models $M_1$ and $M_2$, characterized by the density functions $f_1(x)$ and $f_2(x)$. Let $K_{12} = \int [\log(\frac{f_1(x)}{f_2(x)})] f_1(x) dx$ denote the Kullback-Leibler distance between these two distributions and let $V_{12} = \int [\log(\frac{f_1(x)}{f_2(x)})]^2 f_1(x) dx$. Further assume that $V_{12}$ is finite. We consider the case when no prior information is available and our choice of the model depends only on the likelihood of the data under $M_1$ and $M_2$.

Assume now that a sequence of i.i.d. data $X_1, \ldots, X_m$ is generated according to the model $M_1$. Let $L_1 = \prod_{i=1}^{m} f_1(X_i)$ denote the corresponding likelihood. Let $L_2 = \prod_{i=1}^{m} f_2(X_i)$ denote the likelihood of the data under the wrong model $M_2$. The probability that the likelihood points at the the wrong model $M_2$ is equal to $P(\log L_1 < \log L_2)$. Let us denote $D_{12} = \log L_1 - \log L_2$. Note that

$$D_{12} = \sum_{i=1}^{m} \log \frac{f_1(X_i)}{f_2(X_i)}.$$

Thus, by the Central Limit Theorem, for sufficiently large $m$ the distribution of $D_{12}$ can be approximated by the normal distribution with mean equal to $mK_{12}$ and variance $mV_{12}$.

Consider now the case when the models $M_1$ and $M_2$ belong to the class of mixture densities specified in (2.1). The parameters for the model $M_1$ are $p = 0.01$, $\sigma_1 = 1$ and $\tau_1 = \sqrt{2 \log 200} \approx 3.26$ and the corresponding parameters for the model $M_2$ are $p_2 = 1$, $\sigma_2 = 1$ and $\tau_2 = \sqrt{p}\tau_1 \approx 0.326$. The parameters for the model $M_2$ are chosen in such a way that the probability distributions corresponding to $M_1$ and $M_2$ have the same variance. We used the Monte Carlo method to calculate $K_{12} \approx 0.083$ and $V_{12} \approx 0.33$. Thus $\mathbf{E}(D_{12}) \approx 3.93$, $\mathbf{Var}(D_{12}) \approx 66.58$ and a probability of making a wrong decision $P(D_{12} < 0) \approx 0.31$. Note that while the Kullback Leibler distance between $M_1$ and $M_2$ is rather small these models are completely different in the percentage of alternative hypothesis and the resulting testing procedures give very different results. Wrong decision of accepting the model $M_2$ leads to a rejection of all null hypothesis, while in reality about 99% of them are true. Interestingly, the probability of wrongly detecting the corresponding "full" model $M_2$ quickly decreases with an increase of $p$. This dependence is demonstrated in Figure 5(a). The described phenomenon appears when $\sigma$ is known and unknown and forces us to use the informative prior distribution on $p$ when testing is performed in the sparse mixture setting.

In case when $\sigma$ is unknown we observe a parallel problem with identifying the parameters of the mixture distributions with large values of $p$. For example consider the model $M_1$ with $p = 0.95$, $\sigma_1 = 1$ and $\tau_1 = \sqrt{2 \log 200} \approx 3.26$ and the corresponding "null" model $M_2$ with $p_2 = 0$ and $\sigma_2 = \sqrt{\sigma_1^2 + p\tau_1^2} \approx 3.33$. For this example $K_{12} \approx 0.0013$, $V_{12} \approx 0.0505$ and a probability of making an error $P(D_{12} < 0) \approx 0.37$. Similarly as before, choosing $M_2$ instead of $M_1$ leads to a completely wrong testing procedure (i.e. accepting all hypotheses). In this situation, probability of making a wrong decision increases with $p$ and is illustrated in Figure 5(b). The described phenomenon causes the power of our testing procedures to diminish when the fraction of alternatives exceeds a certain threshold value, as demonstrated in Figure 2(d).

(a) wrong choice of the model with $p = 0$     (b) wrong choice of the model with $p = 0$

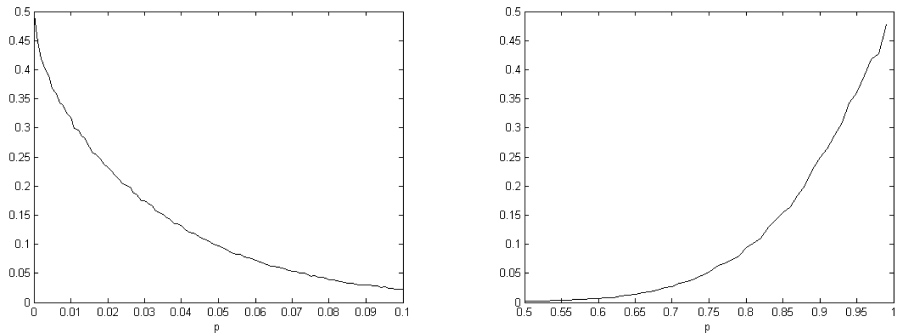

FIG 5. *Probability of making a wrong model choice as a function of true $p$.*

## 7. Conclusions

We have examined several multiple testing procedures keeping in mind both FDR and decision theoretic criteria like MP (Misclassification Probability), efficiency (MP of oracle/MP) and power. We also studied the robustness to some deviations from the assumed prior distribution and compared our fully parametric methods with methods based on nonparametric priors/mixing distributions as in Dirichlet process mixtures or Newton's algorithm.

We observed that if $\sigma$ is known then most methods tend to perform poorly at one of the two extremes. The MLE-based methods (PEB1, BH1) suffer near $p = 0$ due to near lack of identifiability. On the other hand procedures that make use of a conservative prior on $p$ (PEB2, BH2, SB) tend to be too conservative near $p = 1$. Surprisingly, the NPBN procedure, based on Newton's algorithm, does well over the entire range.

Our results confirm the observation of Genovese and Wasserman [12] that for very small $p$'s (for $m = 200$, $p < 0.03$) the misclassification error of BH is close to optimal in the sense of attaining a Bayesian oracle. In this range of $p$ BH works similarly to the Bayes oracle also in terms of FDR and the power. However for $p > 0.03$ the Bayes oracle becomes much more liberal than BH and allows to obtain much smaller misclassification rate. Interestingly, the misclassification probability of the modified version of BH, which uses the knowledge on $p$, is comparable to the misclassification probability of the Bayes oracle over the entire range of $p$. This happens even though the False Discovery Rate and the power of these two are quite different. Our simulations demonstrate that Empirical Bayes methods can be used to estimate $p$ and construct modified versions of BH when the model parameters are unknown.

An interesting methodological fact is that in case when $\sigma$ is unknown all of the considered procedures break down for relatively large $p$. It is somewhat unexpected that one would have a problem when $p$, i.e. the proportion of signals, is large. Section 6 explains how this arises due to the nonidentifiability of mixtures.

The above facts have interesting as well as useful implications for the two applications discussed in Section 2. If our main interest is in QTL and polygenic effects, then $\sigma$ is due to random effects and can be estimated well by appropriate replication. This will virtually reduce the case of unknown $\sigma$ to the known $\sigma$ case and improve the quality of inference. On the other hand, if our goal is QTL mapping alone, then $\sigma$ represents both random effects and polygenic effects and hence can not be directly estimated even with replication. But fortunately for the range of $p$ that is relevant for QTL mapping, namely $p < 0.2$, unknown $\sigma$ does not cause problems at least for $m \geq 200$ (see Figure 2(c)).

## Appendix

**Theorem A.1.** *Let $X_1, \ldots, X_m$ be the sequence of i.i.d. rv's with the density specified by (2.1). Assume that the unknown vector of parameters $\theta_0 = (p_0, \sigma_0^2, \tau_0^2)$ is in the interior of the parameter space $\Omega = [0, 1] \times R^+ \times R^+$. Moreover, assume that the prior density is continuous and positive at $\theta_0$ and that there exists $m_0 \in \mathbf{N}$ such that the corresponding posterior distribution $\Pi(\cdot | x_1, \ldots, x_m)$ is proper when $m \geq m_0$. Then the posterior distribution is consistent, i.e. with probability 1 for* *every Euclidean neighborhood U of $\theta_0$ it holds*

(A.1) $$\{\Pi(U|X_1,\ldots,X_m)\} \to 1 \text{ as } m \to \infty.$$

Theorem 4.4.1 of Ghosh and Ramamoorthi [14] shows that the posterior probability of any weak neighborhood of the mixture distribution $P_{\theta_0}$ tends to one as $m$ tends to infinity. However, the same result for an Euclidean neighborhood of the true parameter $\theta_0$ requires considerably more work. We omit the proof to save space.

**Theorem A.2.** *The misclassification probability of the full Bayes procedure SB, specified in (4.4), converges to the optimal misclassification probability provided by the Bayes oracle (3.9).*

Theorem A.2 essentially follows from Theorem A.1 and regularity properties of the mixture density, but the full proof, though along standard lines, is also somewhat long and hence omitted.

**Theorem A.3.** *The joint posterior distribution of $(\mu_1,\ldots,\mu_m,p,\tau^2,\sigma^2,c)$ under the DPP specification is proper.*

*Proof.* We need to show that

$$I = \int \left[\prod_{i=1}^{m} \frac{1}{\sqrt{2\pi\sigma^2}} \exp\left(-\frac{(x_i-\mu_i)^2}{2\sigma^2}\right)\right] \pi(d\mu_1,\ldots,d\mu_m,d\sigma^2,d\tau^2,dp)$$

is finite. Integrating out $P$ in the hierarchical specification of DPP leads to the following joint conditional distribution of the $\mu_i$'s:

$$\mu_1,\ldots,\mu_m | \sigma^2, \tau^2, p$$
$$\sim \prod_{i=1}^{m} \left[\frac{c}{c+i-1}(1-p)\delta_{\{0\}} + \frac{c}{c+i-1}pN(0,\tau^2) + \frac{1}{c+i-1}\sum_{j=1}^{i-1}\delta_{\{\mu_j\}}\right].$$

From this it can be shown that $I$ equals

$$\sum a(S_0,\Phi) \int_0^\infty \int_0^\infty \prod_{i\in S_0} \frac{1}{\sqrt{2\pi\sigma^2}} \exp\left(-\frac{x_i^2}{2\sigma^2}\right)$$
$$\times \prod_{S\in\Phi} \frac{1}{(\sqrt{2\pi\sigma^2})^{|S|-1}} \exp\left(-\frac{t_S^2}{2\sigma^2}\right) \frac{1}{\sqrt{2\pi(\sigma^2+|S|\tau^2)}}$$
$$\times \exp\left(-\frac{|S|m_S^2}{2(\sigma^2+|S|\tau^2)}\right) \frac{1}{(\sigma^2+\tau^2)^2} d\tau^2 d\sigma^2$$

where the sum is taken over all $S_0 \subset \{1,\ldots,m\}$, $\Phi \in \mathcal{P}(\{1,\ldots,m\}\setminus S_0)$ – the collection of all partitions of $\{1,\ldots,m\}\setminus S_0$,

$$a(S_0,\Phi) = \int \frac{\left[\prod_{i=1}^{|S_0|}(c(1-p)+i-1)\right]\left[\prod_{S\in\Phi}(cp)\{(|S|-1)!\}\right]}{\prod_{i=1}^{m}(c+i-1)} \pi(dp,dc)$$

and

$$m_S = \frac{1}{|S|}\sum_{i\in S} x_i, \quad t_S = \sum_{i\in S}(x_i - m_S)^2.$$

From this one can show that $I < \infty$ by a direct verification of the finiteness of each of the integrals entering the above finite sum. This exercise can be carried out by 1. substituting $Z = \frac{\tau^2}{\sigma^2}$, 2. integrating out $\sigma^2$ and 3. by using the fact that $\tau^2/\sigma^2$ admits a proper density. □

**Acknowledgments.** Much of the work on the manuscript was done during the visits of the first and the third author at Purdue University. The hospitality of the Department of Statistics at Purdue is gratefully appreciated. We also would like to thank the Associate Editor Professor E. Peña for helpful suggestions.